\documentclass[11pt]{amsart}

\usepackage{amsfonts}
\usepackage{amssymb}
\usepackage{hyperref}
\usepackage{stmaryrd}

\usepackage{tikz}
\usetikzlibrary{matrix,arrows,decorations.pathmorphing}

\theoremstyle{plain}
\newtheorem{thm}{Theorem}[section]

\newtheorem{cor}[thm]{Corollary}

\theoremstyle{definition}

\theoremstyle{remark}

\newcommand{\lemref}[1]{\hyperref[#1]{Lemma \ref*{#1}}}
\newcommand{\thmref}[1]{\hyperref[#1]{Theorem \ref*{#1}}}
\newcommand{\propref}[1]{\hyperref[#1]{Proposition \ref*{#1}}}
\newcommand{\corref}[1]{\hyperref[#1]{Corollary \ref*{#1}}}
\newcommand{\defref}[1]{\hyperref[#1]{Definition \ref*{#1}}}
\newcommand{\remref}[1]{\hyperref[#1]{Remark \ref*{#1}}}
\newcommand{\conjref}[1]{\hyperref[#1]{Conjecture \ref*{#1}}}

\newcommand{\Gal}{\mathrm{Gal}}

\def \F {\mathbb{F}}

\makeatletter
\newcommand*{\defeq}{\mathrel{\rlap{%
                     \raisebox{0.27ex}{$\m@th\cdot$}}%
                     \raisebox{-0.27ex}{$\m@th\cdot$}}%
                     =}

\numberwithin{equation}{section}
\makeatother

\makeatletter
\def\@setcopyright{}
\def\serieslogo@{}
\makeatother

\input xy
\xyoption{all} 

\title{Balanced presentations for fundamental groups of curves over finite fields}

\author{Mark Shusterman}
\address{Department of Mathematics, University of Wisconsin-Madison, 480 Lincoln Drive, Madison, WI 53706, USA}
\email{mshusterman@wisc.edu}

\begin{document}

\begin{abstract}

We show that the algebraic fundamental group of a smooth projective curve over a finite field admits a finite topological presentation where the number of relations does not exceed the number of generators.

\end{abstract}

\maketitle

\section{Introduction}

Let $p$ be a prime number, let $k$ be a finite field of characteristic $p$, and let $X$ be a smooth projective curve defined over $k$. 
Associated to $X$ is its algebraic (or \'{e}tale) fundamental group
$
\pi_1(X).
$
This profinite group fits into the short exact sequence
\begin{equation} \label{FES}
1 \to \pi_1(X_{\bar k}) \to \pi_1(X) \to \Gal(\bar k/k) \to 1
\end{equation}
where $\bar k$ is an algebraic closure of $k$, and $X_{\bar k}$ is the base change of $X$ to $\bar k$.

Grothendieck has shown in \cite[exp. XIII]{SGA1} that $\pi_1(X_{\bar k})$, and thus also $\pi_1(X)$, is topologically finitely generated.
Moreover, he shows that for every prime $\ell \neq p$, the maximal pro-$\ell$ quotient of $\pi_1(X_{\bar k})$ admits the pro-$\ell$ presentation
\begin{equation}
\langle x_1, \dots, x_{2g} \mid [x_1, x_2] \cdots [x_{2g-1}, x_{2g}] = 1 \rangle
\end{equation}
where $g$ is the genus of $X$, and $[a,b] = aba^{-1}b^{-1}.$
This result is complemented by the work \cite{Sha47} of Shafarevich, showing that the maximal pro-$p$ quotient of $\pi_1(X_{\bar k})$ is free pro-$p$.

The structure of $\pi_1(X_{\bar k})$, and therefore also of $\pi_1(X)$, is well-understood once $g \leq 1$.
On the other hand, in case $g \geq 2$, the group $\pi_1(X)$ encompasses the structure of $X$ in a non-trivial way. 
One evidence for this is the work \cite{Tam04} of Tamagawa which shows that there are only finitely many smooth projective curves $Y$ over $\bar k$ with $\pi_1(Y) \cong \pi_1(X_{\bar k})$. Another result in this vein, obtained by Mochizuki in \cite{Moch07}, says (roughly speaking) that $X$ can be (functorially) reconstructed from $\pi_1(X)$.

In this work, we shed further light on the structure of $\pi_1(X)$.

\begin{thm} \label{Res}

The \'{e}tale fundamental group of a smooth projective curve over a finite field is topologically finitely presented.
Moreover, the number of relations is at most the number of generators (in some presentation).

\end{thm}

Recall that if $K$ is the function field of $X$, then $\pi_1(X)$ can be identified with $\Gal(K^{\mathrm{ur}}/K)$, where $K^{\mathrm{ur}}$ is the maximal unramified extension of $K$ (inside a fixed separable closure). As a result, the group $\pi_1(X)$ is also studied as a function field analog of the generalized class group $\Gal(L^{\mathrm{ur}}/L)$ of a number field $L$. Examples include the works \cite{BW17, W17} by Boston and Wood. 
These works suggest, vaguely speaking, that for a randomly chosen $X$, the group $\pi_1(X)$ is given by a random balanced presentation (that is, a presentation where the number of generators is at least the number of relations). Motivated by this, Liu and Wood in \cite{LW17} study random balanced presentations and their variants.
\thmref{Res} may thus be viewed as a deterministic counterpart of the works by Boston, Liu, and Wood.

Arithmetic topology, as studied for instance in \cite{Mor11}, postulates that the group $\Gal(L^{\mathrm{ur}}/L)$ should exhibit properties similar to those of the fundamental group of a $3$-manifold. Indeed, balanced presentations are attributes of $3$-manifold groups (see \cite{AFW15}), and Pardon asks in \cite{Par11} whether an analog of this can be established for $\Gal(L^{\mathrm{ur}}/L)$. This serves as an additional motivation for our theorem.

Our proof of \thmref{Res} starts from a formula proved by Lubotzky in \cite{Lub01}. The latter, once combined with Grothendieck's finite generation result mentioned above, reduces our task to an estimation of the dimensions of several cohomology groups. We perform these estimations using the Lyndon-Hochschild-Serre spectral sequence. An important input is duality in cohomology, deduced from the aforementioned results on maximal pro-$\ell$ and pro-$p$ quotients by Grothendieck and Shafarevich.
Consequently, we determine the deficiency of (any presentation of) $\pi_1(X)$, obtaining a result that mirrors Epstein's work \cite{Eps61} on the deficiency of $3$-manifold groups.
Furthermore, our arguments apply not only to $\pi_1(X)$ itself, but to any closed topologically finitely generated subgroup of it.

\begin{cor} \label{cor}

The group $\pi_1(X)$ is topologically coherent.

\end{cor}

Coherence means that all the closed topologically finitely generated subgroups of $\pi_1(X)$ are topologically finitely presented.
Once again, we have an analogy to a result on $3$-manifolds, namely the coherence of their fundamental groups, established by Scott in \cite{Sco73a, Sco73b}.

\section{Proof of \thmref{Res}}

We set $G \defeq \pi_1(X)$, and $N \defeq \pi_1(X_{\bar k})$, so equation \eqref{FES} reads
\begin{equation} \label{RFES}
1 \to N \to G \to  \widehat{ \mathbb{Z}} \to 1.
\end{equation}
Following \cite{Lub01}, denote by $d(G)$ the least cardinality of a generating set of $G$, 
and by $r(G)$ the least number of relations needed to present $G$ (with any number of generators).
Furthermore, define the deficiency of $G$ to be
\begin{equation}
\mathrm{def}(G) \defeq d(G) - r(G).
\end{equation}

If $g = 0$, we get from \cite[Theorem 10.1.2 (i) a]{NSW13} that $N = 1$, so $G \cong \widehat{\mathbb{Z}}$ by \eqref{RFES}. 
Hence $\mathrm{def}(G) = d(G) - r(G) = 1 - 0 = 1$. 
We assume henceforth that $g \geq 1$, and the core of the proof is devoted to showing that $\mathrm{def}(G) = 0$, or equivalently that $r(G) = d(G)$.

By \cite[Theorem 0.2]{Lub01}, we know that $r(G)$ equals
\begin{equation} \label{LF}
\sup_{\lambda} \ \sup_M \Bigg\lceil \frac{\dim_{\F_\lambda} H^2(G,M) - \dim_{\F_\lambda} H^1(G,M)}{\dim_{\F_\lambda} M} \Bigg\rceil + d(G) - \textbf{1}_{M \ncong \F_{\lambda}}
\end{equation}
where $\lambda$ ranges over the prime numbers, $\F_\lambda$ is the finite field with $\lambda$ elements (and a trivial $G$-action),
and $M$ ranges over all simple $\mathbb{F}_\lambda \llbracket G \rrbracket$-modules.

Let us first see that $r(G) \geq d(G)$. For that, take any $\lambda \neq p$ and let $M$ be the (trivial) $G$-module $\F_\lambda$.
Grothendieck's aforementioned result on the maximal pro-$\lambda$ quotient of (open subgroups of) $N$
implies that $N$ is a Poincar\'{e} group of dimension $2$ at $\lambda$ 
(for an alternative proof, see \cite[Theorem 10.1.2 (i) b]{NSW13}). As $\widehat{ \mathbb{Z}}$ is a Poincar\'{e} group of dimension $1$ at $\lambda$, it follows from \cite[Theorem 3.7.4]{NSW13} applied to \eqref{RFES} that $G$ is a Poincar\'{e} group of dimension $3$ at $\lambda$ (see also \cite[Corollary 10.1.3 (ii)]{NSW13}). 
The latter implies that
\begin{equation}
\dim_{\F_\lambda} H^2 \big(G,M\big) = 
\dim_{\F_\lambda} H^{1}\big(G, \mathrm{Hom}(M, \F_\lambda)\big) = 
\dim_{\F_\lambda} H^1\big(G,M\big)
\end{equation}
as $M = \F_\lambda$. It follows from \eqref{LF} that indeed $r(G) \geq d(G)$.

We shall now start proving that $r(G) \leq d(G)$. For that, we fix a prime $\lambda$ and a simple (and thus finite) $\F_\lambda \llbracket G \rrbracket$-module $M$. 
By \eqref{LF}, it suffices to show that
\begin{equation} \label{Goal}
\dim_{\F_\lambda} H^2(G,M) - \dim_{\F_\lambda} H^1(G,M) \leq \textbf{1}_{M \ncong \F_{\lambda}} \cdot \dim_{\F_\lambda} M.
\end{equation}

The Lyndon-Hochschild-Serre spectral sequence, associated to \eqref{RFES} and to $M$, gives rise (see \cite{KHW12}) to the seven-term exact sequence
\begin{equation*}
\begin{split}
0 &\to H^1\big(\widehat{\mathbb{Z}}, M^N\big) \to H^1\big(G,M\big) \to H^1\big(N,M\big)^{\widehat{\mathbb{Z}}} \to 
H^2\big(\widehat{\mathbb{Z}}, M^N\big) \\
&\to \mathrm{Ker}\big(H^2(G,M) \rightarrow H^2(N,M) \big) \to H^1\big(\widehat{\mathbb{Z}},H^1(N,M)\big) \to H^3\big(\widehat{\mathbb{Z}},M^N\big).
\end{split}
\end{equation*}
Since the cohomological dimension of $\widehat{\mathbb{Z}}$ is $1$, 
the rightmost term in each row of the exact sequence above vanishes, so we obtain the exact sequence
\begin{equation} \label{H1}
0 \to H^1\big(\widehat{\mathbb{Z}}, M^N\big) \to H^1\big(G,M\big) \to H^1\big(N,M\big)^{\widehat{\mathbb{Z}}} \to 0
\end{equation}
and the exact sequence
\begin{equation} \label{H2}
0 \to  H^1\big(\widehat{\mathbb{Z}},H^1(N,M)\big) \to H^2\big(G, M\big) \to H^2\big(N, M\big).
\end{equation}


It follows from \eqref{H1} that
\begin{equation} \label{h1}
\dim_{\F_\lambda} H^1\big(G,M\big) \geq \dim_{\F_\lambda}  H^1\big(N,M\big)^{\widehat{\mathbb{Z}}}
\end{equation}
and from \eqref{H2} (combined with \cite[Proposition 1.7.7 (i)]{NSW13}) that
\begin{equation} \label{h2}
\begin{split}
\dim_{\F_\lambda} H^2\big(G,M\big) &\leq 
\dim_{\F_\lambda}  H^1\big(\widehat{\mathbb{Z}},H^1(N,M)\big) + \dim_{\F_\lambda} H^2\big(N, M\big) \\
&\leq \dim_{\F_\lambda} H^1\big(N,M\big)_{\widehat{\mathbb{Z}}} + \dim_{\F_\lambda} H^2\big(N, M\big).
\end{split}
\end{equation}
In order to establish \eqref{Goal}, we want to compare \eqref{h1} with \eqref{h2}, so we now show that
\begin{equation} \label{TP}
\dim_{\F_\lambda}  H^1\big(N,M\big)^{\widehat{\mathbb{Z}}} = \dim_{\F_\lambda} H^1\big(N,M\big)_{\widehat{\mathbb{Z}}}.
\end{equation}

From Grothendieck's result, we know that $N$ is finitely generated, so $V \defeq H^1(N,M)$ is a finite dimensional vector space over $\F_\lambda$. Therefore, taking a generator $\sigma$ of $\widehat{\mathbb{Z}}$, and considering the linear map it induces on $V$, we find that
\begin{equation} \label{LA}
\begin{split}
\dim_{\F_\lambda} V_{\widehat{\mathbb{Z}}}  &= \dim_{\F_\lambda} V/\mathrm{Im}(\sigma - 1)
= \dim_{\F_\lambda} V - \dim_{\F_\lambda} \mathrm{Im}(\sigma - 1) \\
&= \dim_{\F_\lambda} \mathrm{Ker}(\sigma - 1) = \dim_{\F_\lambda} V^{\widehat{\mathbb{Z}}}
\end{split}
\end{equation}
so we have arrived at \eqref{TP}. 

Subtracting \eqref{h1} from \eqref{h2}, we can therefore use \eqref{TP} to get that
\begin{equation} \label{BD}
\dim_{\F_\lambda} H^2\big(G,M\big) - \dim_{\F_\lambda} H^1\big(G,M\big) \leq 
\dim_{\F_\lambda} H^2\big(N,M\big).
\end{equation}
We are now ready to verify \eqref{Goal}, distinguishing between several cases. 

Suppose that $\lambda \neq p$. We have already treated the case $M = \F_\lambda$ while showing that $r(G) \geq d(G)$, so let us assume that $\textbf{1}_{M \ncong \F_\lambda} = 1$. As already mentioned before, $N$ is a Poincar\'{e} group of dimension $2$ at $\lambda$, so we have
\begin{equation} \label{CD}
\begin{split}
\dim_{\F_\lambda} H^2\big(N,M\big) &= \dim_{\F_\lambda} H^0\big(N, \mathrm{Hom}(M, \F_\lambda) \big) \\
&\leq \dim_{\F_\lambda} \mathrm{Hom}\big(M, \F_\lambda\big) = \dim_{\F_\lambda} M.
\end{split}
\end{equation}
Hence, in this case, \eqref{Goal} is a consequence of \eqref{BD} and \eqref{CD}.

Suppose now that $\lambda = p$.
The freeness of the maximal pro-$\lambda$ quotient of (open subgroups of) $N$ proved by Shafarevich,
implies that the cohomological $\lambda$-dimension of $N$ is $1$ (for an alternative proof, see \cite[Theorem 10.1.12 (ii)-(iii)]{NSW13}).
This implies that
\begin{equation} \label{DD}
\dim_{\F_\lambda} H^2\big(N,M\big) = 0.
\end{equation}
Therefore, in this case, \eqref{Goal} follows from \eqref{BD} and \eqref{DD}.
We have thus finished showing that $r(G) \leq d(G)$. 

By \cite[Corollary 2.5]{Lub01}, there exists a presentation of $G$ with $d(G)$ generators and $r(G)$ relations.
We have shown that $r(G) = d(G)$, so this is a balanced presentation. 
We have seen that $N$ is not procyclic, therefore
\begin{equation}
r(G) = d(G) > 1
\end{equation}
so in the language of \cite[Theorem 0.2]{Lub01}, $G$ is not $d(G)$-abelian-indexed. 
Hence, it follows from \cite[Theorem 0.1]{Lub01} that the deficiency of every presentation of $G$ equals $\mathrm{def}(G)$.
That is, for every presentation
\begin{equation}
1 \to K \to F \to G \to 1
\end{equation}
where $F$ is a free profinite group, we have
\begin{equation}
d(F) - d_F(K) = \mathrm{def}(G) = 0.
\end{equation}
Here, $d_F(K)$ is the number of relations in the presentation, namely the least cardinality of a generating set for $K$ as a normal subgroup of $F$.

\section{Deducing \corref{cor}}

An extension of the above proof that $r(G) \leq d(G)$ shows that any finitely generated $G_0 \leq_c G$ is finitely presented. 
Indeed, setting $N_0 \defeq G_0 \cap N$ we get the exact sequence
\begin{equation}
1 \to N_0 \to G_0 \to C \to 1
\end{equation}
where $C$ is a procyclic group.
Fixing a prime number $\lambda$, and a simple $\F_\lambda \llbracket G_0 \rrbracket$-module $M$, 
we arrive (arguing as in \eqref{h1} and \eqref{h2}) at the bound
\begin{equation*}
\dim H^2\big(G_0,M\big) - \dim H^1\big(G_0,M\big) \leq 
\dim H^2\big(N_0,M\big) + \dim V_{C} - \dim V^{C}
\end{equation*}
where (as before) $V \defeq H^1(N_0, M)$, and dimensions are taken over $\F_\lambda$.

Let $W$ be a finite dimensional quotient of $V$ on which $C$ acts trivially. 
As $V$ is a direct limit of finite $C$-modules, there exists a finite dimensional $C$-submodule $U$ of $V$ that surjects onto $W$.
Arguing as in \eqref{LA}, we get
\begin{equation}
\dim W \leq \dim U_C = \dim U^C \leq \dim V^C
\end{equation}
so it follows that 
\begin{equation} \label{WE}
\dim V_C \leq \dim V^C.
\end{equation}
Since $G_0$ is finitely generated, it follows from (the analog of) \eqref{h1} that $\dim V^C$ is finite,
so combining the first inequality of this section with \eqref{WE}, we obtain
\begin{equation} \label{BDD}
\dim H^2\big(G_0,M\big) - \dim H^1\big(G_0,M\big) \leq \dim H^2\big(N_0,M\big)
\end{equation}
which is the analog of \eqref{BD}.

In \eqref{CD} and \eqref{DD}, it is shown that
\begin{equation} \label{H2BB}
\dim H^2(L, M) \leq \dim M
\end{equation}
for every open subgroup $L$ of $N$ (to which the action of $N_0$ on $M$ extends).
Writing $N_0$ as the inverse limit of a family of open subgroups of $N$, 
we conclude from \eqref{H2BB} (using \cite[Proposition 1.5.1]{NSW13}) that
\begin{equation} \label{H2BBB}
\dim H^2(N_0, M) \leq \dim M.
\end{equation}
Combining \eqref{H2BBB} with \eqref{BDD}, we conclude from \eqref{LF} (or from \cite[Theorem 0.3]{Lub01}) that $G_0$ is finitely presented.

\section*{Acknowledgments}

Mark Shusterman is grateful to Melanie Matchett Wood for all the discussions and encouragement in the course of the work on this paper.

\end{document}